\documentclass{amsart}
\usepackage{amsmath,amssymb}
\usepackage{yhmath}

\newcommand{\bdis}{\begin{displaymath}}
\newcommand{\edis}{\end{displaymath}}
\newcommand{\be}{\begin{equation}}
\newcommand{\ee}{\end{equation}}
\newcommand{\mbb}{\mathbb}
\newcommand{\mcal}{\mathcal} 
\newcommand{\mfrak}{\mathfrak}

\newcommand{\vp}{\varphi}

\newcommand{\zf}{\zeta\left(\frac{1}{2}+it\right)}

\theoremstyle{definition}

\theoremstyle{remark}
\newtheorem{remark}[]{Remark}

\newtheorem*{mydef11}{{\bf Theorem 1}}

\newtheorem*{mydef12}{{\bf Theorem 2}}

\newtheorem*{mydef13}{{\bf Theorem 3}}

\newtheorem*{mydef14}{{\bf Theorem 4}}

\newtheorem*{mydef15}{{\bf Theorem 5}} 

\newtheorem*{mydef16}{{\bf Theorem 6}} 

\newtheorem*{mydef17}{{\bf Theorem 7}}

\newtheorem*{mydef51}{{\bf Lemma 1}}

\newtheorem*{mydef52}{{\bf Lemma 2}}

\newtheorem*{mydef53}{{\bf Lemma 3}}

\newtheorem*{mydef54}{{\bf Lemma 4}}

\newtheorem*{mydef55}{{\bf Lemma 5}}

\newtheorem*{mydef56}{{\bf Lemma 6}}

\newtheorem*{mydef57}{{\bf Lemma 7}} 

\newtheorem*{mydef58}{{\bf Lemma 8}} 

\newtheorem*{mydef59}{{\bf Lemma 9}} 

\newtheorem*{mydef510}{{\bf Lemma 10}} 

\newtheorem*{mydef511}{{\bf Lemma 11}} 

\newtheorem*{mydef512}{{\bf Lemma 12}} 

\newtheorem*{mydef513}{{\bf Lemma 13}}

\newtheorem*{mydef81}{{\bf Property 1}}

\newtheorem*{mydef82}{{\bf Property 2}}

\newtheorem*{mydefMR}{{\bf Main result}}

\numberwithin{equation}{section}

\begin{document}

\title[Jacob's ladders and three new equivalents \dots]{Jacob's ladders and three new equivalents of the Fermat-Wiles theorem and an infinite set of these equivalents that are independent on the Jacob's ladders}

\author{Jan Moser}

\address{Department of Mathematical Analysis and Numerical Mathematics, Comenius University, Mlynska Dolina M105, 842 48 Bratislava, SLOVAKIA}

\email{jan.mozer@fmph.uniba.sk}


\begin{abstract}
In this paper we show that there is an infinite set of points of contact between the set of all Dirichlet's series and Fermat-Wiles theorem. The proof is independent on the Jacob's ladders. 
\end{abstract}
\maketitle

\section{Introduction} 

\subsection{} 

In this continuation of the series of papers \cite{8} -- \cite{14} we use our method to construct new functionals on: 

\begin{itemize}
	\item[(A)] the basic integral 
	\be \label{1.1} 
	\int_1^T|\zeta(\sigma+it)|^2{\rm d}t\sim \zeta(2\sigma)T,\ \sigma>1,\ T\to \infty , 
	\ee 
	that is connected directly with the definition of the Riemann's zeta-function in $\mbb{C}$ (1859): 
	\be \label{1.2} 
	\zeta(s)=\sum_{n=1}^\infty \frac{1}{n^s}=\prod_p\frac{1}{1-\frac{1}{p^s}},\ s=\sigma+it,\ \sigma>1; 
	\ee 
	\item[(B)] the classical Hardy-Littlewood continuation (1922) of the formula (\ref{1.1}) into the critical strip 
	\be \label{1.3} 
	\frac 12<\sigma\leq 1, 
	\ee 
	i. e. onto the region $\sigma>\frac 12$, comp. \cite{16}, pp. 29 -- 31; 
	\item[(C)] the classical Selberg's integral (1946): 
	\be \label{1.4} 
	\int_0^T|S_1(t)|^{2l}{\rm d}t\sim \bar{c}(l)T,\ T\to\infty 
	\ee 
	(see \cite{15}, p. 130) for every fixed $l\in\mbb{N}$, where 
	\be \label{1.5} 
	S_1(t)=\frac{1}{\pi}\int_0^t\arg\zeta\left(\frac 12+iu\right){\rm d}u. 
	\ee 
\end{itemize} 

\subsection{} 

In this paper we have obtained the first proof of the fact, that the formula\footnote{Comp. (\ref{1.1}).} (well known for more than 100 years)
\be \label{1.6} 
\lim_{T\to\infty}\frac{1}{T}\int_0^T|\zeta(\sigma+it)|^2{\rm d}t=\zeta(2\sigma),\ \sigma>1
\ee 
is followed by the existence of the point of contact between the Riemann's zeta-function and the Fermat-Wiles theorem. Our proof is independent on the Jacob's ladders. 

Namely, by means of the formula (\ref{1.6}), we have obtained the following 
\be \label{1.7} 
\lim_{\tau\to\infty}\frac{1}{\tau}\int_1^{\frac{x}{\zeta(2\sigma)}\tau}|\zeta(\sigma+it)|^2{\rm d}t=x 
\ee 
for every fixed $x$ and every fixed $\sigma>1$. 

Now, in the special case of the Fermat's rationals 
\be \label{1.8} 
x\to \frac{x^n+y^n}{z^n},\ x,y,z,n\in\mbb{N},\ n\geq 3, 
\ee 
the condition: 
\be \label{1.9} 
\lim_{\tau\to\infty}\frac{1}{\tau}\int_1^{\frac{x^n+y^n}{z^n}\frac{\tau}{\zeta(2\sigma)}}|\zeta(\sigma+it)|^2{\rm d}t\not=1
\ee 
follows on the class of all Fermat's rationals and this represents new equivalent to the Fermat-Wiles theorem.  

Next, we have obtained similar result also for the classical Selberg's formula\footnote{Comp. (\ref{1.4}).} 
\be \label{1.10} 
\lim_{T\to\infty}\frac{1}{T}\int_0^T|S_1(t)|^{2l}{\rm d}t=\bar{c}(l). 
\ee 
Namely the functional 
\be \label{1.11} 
\lim_{\tau\to\infty}\frac{1}{\tau}\int_0^{\frac{x}{\bar{c}(l)}\tau}|S_1(t)|^{2l}{\rm d}t=x,\ x>0  
\ee 
and corresponding equivalent to the Fermat-Wiles theorem 
\be \label{1.12} 
\lim_{\tau\to\infty}\frac{1}{\tau}\int_0^{\frac{x^n+y^n}{z^n}\frac{\tau}{\bar{c}(l)}}|S_1(t)|^{2l}{\rm d}t\not=1.  
\ee 

\subsection{} 

Here we give some remarks. 

\begin{remark}
We have obtained the conditions (\ref{1.9}) and (\ref{1.12}) independently on the Jacob's ladders. 
\end{remark} 

\begin{remark} 
Formulae (\ref{1.7}), (\ref{1.9}), (\ref{1.11}) and (\ref{1.12}) connected with quite simple results (\ref{1.6}) and (\ref{1.10})  have not engage the attention of Hardy nor Littlewood and Selberg. This can be viewed as missed chance in the past. 
\end{remark} 

\begin{remark}
Our remark 2 is good illustration of the following L. D. Kudryavtsev's (see \cite{2}, Preface) statement: "Formulas can be more reasonable than the ones who use them and can give more than expected." 
\end{remark}

\begin{remark}
Even if the equivalents (\ref{1.9}) and (\ref{1.12}) are independent on the Jacob's ladders, the integrals are tightly coupled with the Jacob's ladders: 
\be \label{1.13} 
\int_1^T|\zeta(\sigma+it)|^2{\rm d}t\sim \int_{\frac{\zeta(2\sigma)}{1-c}T}^{[\frac{\zeta(2\sigma)}{1-c}T]^1}\left|\zf\right|^2{\rm d}t,\ T\to\infty, 
\ee  
\be \label{1.14} 
\int_0^T|S_1(t)|^{2l}{\rm d}t\sim \int_{\frac{\bar{c}(l)}{1-c}T}^{[\frac{\bar{c}(l)}{1-c}T]^1}\left|\zf\right|^2{\rm d}t,\ T\to\infty, 
\ee  
where 
\be \label{1.15} 
\left[\frac{\zeta(2\sigma)}{1-c}T\right]^1=\vp_1^{-1}\left(\frac{\zeta(2\sigma)}{1-c}T\right). 
\ee 
\end{remark} 

\begin{remark}
We give, immediately, also one generalization of the condition (\ref{1.9}) in the Hardy-Littlewood case $\sigma>\frac 12$: 
\be \label{1.16} 
\lim_{\tau\to\infty}\frac{1}{\tau}\int_1^{\frac{x^n+y^n}{z^n}\frac{\tau}{\zeta(2\frac{\bar{x}^m+\bar{y}^m}{\bar{z}^m})}}
\left|\zeta\left(2\frac{\bar{x}^m+\bar{y}^m}{\bar{z}^m}+it\right)\right|^2{\rm d}t\not=1, 
\ee  
where 
\be \label{1.17} 
\frac{\bar{x}^m+\bar{y}^m}{\bar{z}^m}\geq \frac 12+\epsilon;\ \bar{x}, \bar{y}, \bar{z}, m\in\mbb{N},\ m\geq 3, 
\ee 
i. e. we have the formula for a pair of Fermat's rationals 
\be \label{1.18} 
\left(\frac{x^n+y^n}{z^n}, \frac{\bar{x}^m+\bar{y}^m}{\bar{z}^m}\right)
\ee 
with the condition (\ref{1.17}). 
\end{remark} 

\subsection{} 

In this subsection we give one elements of the set of new equivalents of the Fermat-Wiles theorem that is constructed by means of Jacob's ladders. 

First, we obtain the following functional 
\be \label{1.19} 
\begin{split}
& \lim_{\tau\to\infty}\frac{1}{\tau}\int_{\frac{x}{2k(1-c)\bar{c}(l)\zeta(2\sigma)}}^{[\frac{x}{2k(1-c)\bar{c}(l)\zeta(2\sigma)}]^k}
\left\{ 
2\bar{c}(l)\zeta(2\sigma)\left|\zf\right|^2+ \right. \\ 
& \left.(1-c)\bar{c}(l)|\zeta(\sigma+it)|^2+(1-c)\zeta(2\sigma)|S_1(t)|^{2l}
\right\}{\rm d}t=x
\end{split}
\ee  
for every fixed $x>0,\ k,l\in\mbb{N}$, and $\epsilon$ is a small positive number, where 
\be \label{1.20} 
[Y]^k=\vp_1^{-k}(Y). 
\ee  

\begin{remark}
The functional (\ref{1.19}) represents the result of cooperation of the functions 
\be \label{1.21} 
\left|\zf\right|^2,\ |\zeta(\sigma+it)|^2,\ |S_1(t)|^{2l}
\ee 
generated by the Riemann's function $\zeta(s)$, comp. (\ref{1.5}). 
\end{remark} 

Now, it follows from the functional (\ref{1.19}) in the special case of the Fermat's rationals: the condition 
\be \label{1.22} 
\begin{split}
	& \lim_{\tau\to\infty}\frac{1}{\tau}\int_{\frac{x^n+y^n}{z^n}\frac{\tau}{2k(1-c)\bar{c}(l)\zeta(2\sigma)}}^{[\frac{x^n+y^n}{z^n}\frac{\tau}{2k(1-c)\bar{c}(l)\zeta(2\sigma)}]^k}
	\left\{ 
	2\bar{c}(l)\zeta(2\sigma)\left|\zf\right|^2+ \right. \\ 
	& \left.(1-c)\bar{c}(l)|\zeta(\sigma+it)|^2+(1-c)\zeta(2\sigma)|S_1(t)|^{2l}
	\right\}{\rm d}t\not= 1
\end{split}
\ee 
on the class of all Fermat's rationals represents the new equivalent of the Fermat-Wiles theorem. 

\subsection{} 

Immediate generalizations of the results (\ref{1.7}) and (\ref{1.9}) are presented in this subsection. 

Let us reminf the old mean-value theorem for Dirichlet's series: if 
\be \label{1.23} 
f(s)=\sum_{n=1}^\infty\frac{a_n}{n^s},\ s=\sigma+it,\ a_n\in\mbb{C}
\ee 
is absolutely convergent for $\sigma=\sigma_0$, then, for $\sigma=\sigma_0$, 
\be \label{1.24} 
\lim_{T\to\infty}\frac{1}{T}\int_0^T|f(\sigma_0+it)|^2{\rm d}t=\sum_{n=1}^\infty\frac{|a_n|^2}{n^{2\sigma_0}}.  
\ee  

\begin{remark}
Of course, the Riemann's zeta-function is a canonical element of the infinite set of all Dirichlet's series $\mfrak{D}$. 
\end{remark} 

Let us put 
\be \label{1.25} 
\sum_{n=1}^\infty\frac{|a_n|^2}{n^{2\sigma_0}}=F(\sigma_0;f), 
\ee  
where 
\be \label{1.26} 
F(\sigma_0;f)=F[\sigma_0(f);f]\in\mbb{R}^+
\ee 
for every fixed element $f$ in $\mfrak{D}$ and corresponding point of absolute convergence $\sigma_0(f)$. 

Now, it is clear, that the fact 
\be \label{1.27} 
\{(\ref{1.6})\}\ \Rightarrow \ \{(\ref{1.7})\}\wedge \{(\ref{1.9})\}
\ee 
gives the result 
\be \label{1.28} 
\lim_{\tau\to\infty}\frac{1}{\tau}\int_0^{\frac{x}{F(\sigma_0;f)}\tau}|f(\sigma_0+it)|^2{\rm d}t=x,\ x>0
\ee 
and also the following.  

\begin{mydefMR}
The $\mfrak{D}$-condition 
\be \label{1.29} 
\lim_{\tau\to\infty}\frac{1}{\tau}\int_0^{\frac{x^n+y^n}{z^n}\frac{\tau}{F(\sigma_0;f)}}|f(\sigma_0+it)|^2{\rm d}t\not=1
\ee 
on the class of all Fermat's rationals represents new $\mfrak{D}$-equivalent of the Fermat-Wiles theorem for every fixed $f\in\mfrak{D}$ and corresponding $\sigma_0(f)$. 
\end{mydefMR} 

\begin{remark}
In this way we have proved the existence of infinite set of points od contact between the set $\mfrak{D}$ and the Fermat-Wiles theorem. 
\end{remark}

\section{Jacob's ladders: notions and basic geometrical properties}  

\subsection{} 

In this paper we use the following notions of our works \cite{3} -- \cite{7}: 
\begin{itemize}
\item[{\tt (a)}] Jacob's ladder $\vp_1(T)$, 
\item[{\tt (b)}] direct iterations of Jacob's ladders 
\bdis 
\begin{split}
	& \vp_1^0(t)=t,\ \vp_1^1(t)=\vp_1(t),\ \vp_1^2(t)=\vp_1(\vp_1(t)),\dots , \\ 
	& \vp_1^k(t)=\vp_1(\vp_1^{k-1}(t))
\end{split}
\edis 
for every fixed natural number $k$, 
\item[{\tt (c)}] reverse iterations of Jacob's ladders 
\be \label{2.1}  
\begin{split}
	& \vp_1^{-1}(T)=\overset{1}{T},\ \vp_1^{-2}(T)=\vp_1^{-1}(\overset{1}{T})=\overset{2}{T},\dots, \\ 
	& \vp_1^{-r}(T)=\vp_1^{-1}(\overset{r-1}{T})=\overset{r}{T},\ r=1,\dots,k, 
\end{split} 
\ee   
where, for example, 
\be \label{2.2} 
\vp_1(\overset{r}{T})=\overset{r-1}{T}
\ee  
for every fixed $k\in\mbb{N}$ and every sufficiently big $T>0$. We also use the properties of the reverse iterations listed below.  
\be \label{2.3}
\overset{r}{T}-\overset{r-1}{T}\sim(1-c)\pi(\overset{r}{T});\ \pi(\overset{r}{T})\sim\frac{\overset{r}{T}}{\ln \overset{r}{T}},\ r=1,\dots,k,\ T\to\infty,  
\ee 
\be \label{2.4} 
\overset{0}{T}=T<\overset{1}{T}(T)<\overset{2}{T}(T)<\dots<\overset{k}{T}(T), 
\ee 
and 
\be \label{2.5} 
T\sim \overset{1}{T}\sim \overset{2}{T}\sim \dots\sim \overset{k}{T},\ T\to\infty.   
\ee  
\end{itemize} 

\begin{remark}
	The asymptotic behaviour of the points 
	\bdis 
	\{T,\overset{1}{T},\dots,\overset{k}{T}\}
	\edis  
	is as follows: at $T\to\infty$ these points recede unboundedly each from other and all together are receding to infinity. Hence, the set of these points behaves at $T\to\infty$ as one-dimensional Friedmann-Hubble expanding Universe. 
\end{remark}  

\subsection{} 

Let us remind that we have proved\footnote{See \cite{7}, (3.4).} the existence of almost linear increments 
\be \label{2.6} 
\begin{split}
& \int_{\overset{r-1}{T}}^{\overset{r}{T}}\left|\zf\right|^2{\rm d}t\sim (1-c)\overset{r-1}{T}, \\ 
& r=1,\dots,k,\ T\to\infty,\ \overset{r}{T}=\overset{r}{T}(T)=\vp_1^{-r}(T)
\end{split} 
\ee 
for the Hardy-Littlewood integral (1918) 
\be \label{2.7} 
J(T)=\int_0^T\left|\zf\right|^2{\rm d}t. 
\ee  

For completeness, we give here some basic geometrical properties related to Jacob's ladders. These are generated by the sequence 
\be \label{2.8} 
T\to \left\{\overset{r}{T}(T)\right\}_{r=1}^k
\ee 
of reverse iterations of of the Jacob's ladders for every sufficiently big $T>0$ and every fixed $k\in\mbb{N}$. 

\begin{mydef81}
The sequence (\ref{2.8}) defines a partition of the segment $[T,\overset{k}{T}]$ as follows 
\be \label{2.9} 
|[T,\overset{k}{T}]|=\sum_{r=1}^k|[\overset{r-1}{T},\overset{r}{T}]|
\ee 
on the asymptotically equidistant parts 
\be \label{2.10} 
\begin{split}
& \overset{r}{T}-\overset{r-1}{T}\sim \overset{r+1}{T}-\overset{r}{T}, \\ 
& r=1,\dots,k-1,\ T\to\infty. 
\end{split}
\ee 
\end{mydef81} 

\begin{mydef82}
Simultaneously with the Property 1, the sequence (\ref{2.8}) defines the partition of the integral 
\be \label{2.11} 
\int_T^{\overset{k}{T}}\left|\zf\right|^2{\rm d}t
\ee 
into the parts 
\be \label{2.12} 
\int_T^{\overset{k}{T}}\left|\zf\right|^2{\rm d}t=\sum_{r=1}^k\int_{\overset{r-1}{T}}^{\overset{r}{T}}\left|\zf\right|^2{\rm d}t, 
\ee 
that are asymptotically equal 
\be \label{2.13} 
\int_{\overset{r-1}{T}}^{\overset{r}{T}}\left|\zf\right|^2{\rm d}t\sim \int_{\overset{r}{T}}^{\overset{r+1}{T}}\left|\zf\right|^2{\rm d}t,\ T\to\infty. 
\ee 
\end{mydef82} 

It is clear, that (\ref{2.10}) follows from (\ref{2.3}) and (\ref{2.5}) since 
\be \label{2.14} 
\overset{r}{T}-\overset{r-1}{T}\sim (1-c)\frac{\overset{r}{T}}{\ln \overset{r}{T}}\sim (1-c)\frac{T}{\ln T},\ r=1,\dots,k, 
\ee  
while our eq. (\ref{2.13}) follows from (\ref{2.6}) and (\ref{2.5}). 

\section{The first set of equivalents of the Fermat-Wiles theorem without Jacob's ladders} 

\subsection{} 

We use the basic formula (\ref{1.1}) in the following form 
\be \label{3.1} 
\int_1^T|\zeta(\sigma+it)|^2{\rm d}t=\zeta(2\sigma)T+o(T),\ \sigma>1,\ T\to\infty. 
\ee 
Now, if we put 
\be \label{3.2} 
T=\frac{x}{\zeta(2\sigma)}\tau,\ x>0;\ \{T\to\infty\} \ \Leftrightarrow \ \{\tau\to\infty\}
\ee 
into (\ref{3.1}), then we obtain the following. 

\begin{mydef51}
\be \label{3.3} 
\lim_{\tau\to\infty}\frac{1}{\tau}\int_1^{\frac{x}{\zeta(2\sigma)}\tau}\zeta(\sigma+it)|^2{\rm d}t=x 
\ee  
for every fixed $x>0$ and every fixed $\sigma>1$. 
\end{mydef51} 

\begin{remark}
Of course, the functional 
\be \label{3.4} 
y(\tau;x)=\frac{x}{\zeta(2\sigma)}\tau\to x 
\ee 
is defined by the eq. (\ref{3.3}), where the ray 
\be \label{3.5} 
y(\tau;x)=\frac{x}{\zeta(2\sigma)}\tau,\ \tau\in(\tau_1(x,\sigma),+\infty)
\ee 
is defined by its slope 
\be \label{3.6} 
0<\arctan\frac{x}{\zeta(2\sigma)}<\frac{\pi}{2}, 
\ee  
and $\tau_1(x;\tau)>0$ is sufficiently big for every fixed $x>0$ and every fixed $\sigma>1$. 
\end{remark} 

If we use the substitution\footnote{See (\ref{1.8})} 
\be \label{3.7} 
x\to \frac{x^n+y^n}{z^n}
\ee 
in (\ref{3.3}), then we obtain the following. 

\begin{mydef52}
\be \label{3.8} 
\lim_{\tau\to\infty}\frac{1}{\tau}\int_1^{\frac{x^n+y^n}{z^n}\frac{\tau}{\zeta(2\sigma)}}\zeta(\sigma+it)|^2{\rm d}t=\frac{x^n+y^n}{z^n}
\ee 
for every fixed Fermat's rational and every fixed $\sigma>1$. 
\end{mydef52} 

Consequently, we have the following. 

\begin{mydef11}
The $\zeta$-condition 
\be \label{3.9}  
\lim_{\tau\to\infty}\frac{1}{\tau}\int_1^{\frac{x^n+y^n}{z^n}\frac{\tau}{\zeta(2\sigma)}}|\zeta(\sigma+it)|^2{\rm d}t\not=1 
\ee 
on the class of all Fermat's rationals represents the $\zeta$-equivalent of the Fermat-Wiles theorem for every fixed $\sigma>1$. 
\end{mydef11} 

\subsection{} 

Let us remind the following two Selberg's formulas, see \cite{15}, p. 130: 
\be \label{3.10} 
\int_T^{T+H}|S_1(t)|^{2l}{\rm d}t=\bar{c}(l)H+\mcal{O}\left(\frac{H}{\ln T}\right), 
\ee 
\be \label{3.11} 
\int_0^T|S_1(t)|^{2l}{\rm d}t=\bar{c}{l}T+\mcal{O}\left(\frac{T}{\ln T}\right), 
\ee 
where 
\be \label{3.12} 
T^a\leq H\leq T,\ \frac{1}{2}<a\leq 1,\ \bar{c}(l)>0,\ l\in\mbb{N}. 
\ee 

We use the formula (\ref{3.11}) in this subsection. We put 
\be \label{3.13} 
T=\frac{x}{\bar{c}(l)}\tau,\ x>0 
\ee  
into (\ref{3.11}) to obtain the following. 

\begin{mydef53}
\be \label{3.14} 
\lim_{\tau\to\infty}\frac{1}{\tau}\int_0^{\frac{x}{\bar{c}(l)}\tau}|S_1(t)|^{2l}{\rm d}t=x 
\ee 
for every fixed $x>0$ and $l\in\mbb{N}$.  
\end{mydef53} 

Next, by making use of the substitution (\ref{3.7}) in (\ref{3.14}), we get 

\begin{mydef54}
\be \label{3.15} 
\lim_{\tau\to\infty}\frac{1}{\tau}\int_0^{\frac{x^n+y^n}{z^n}\frac{\tau}{\bar{c}(l)}}|S_1(t)|^{2l}{\rm d}t=\frac{x^n+y^n}{z^n}
\ee 
for every fixed Fermat's rational. 
\end{mydef54} 

Consequently we have the following result. 

\begin{mydef12}
The $\zeta$-condition 
\be \label{3.16} 
\lim_{\tau\to\infty}\frac{1}{\tau}\int_0^{\frac{x^n+y^n}{z^n}\frac{\tau}{\bar{c}(l)}}|S_1(t)|^{2l}{\rm d}t\not=1 
\ee  
on the class of all Fermat's rationals represents the next $\zeta$-equivalent of the Fermat-Wiles theorem for every fixed $l\in\mbb{N}$, where 
\bdis 
S_1(t)=\frac{1}{\pi}\int_0^t\arg\zeta\left( \frac 12+iu\right){\rm d}u. 
\edis 
\end{mydef12} 

\subsection{} 

Since 
\be \label{3.17} 
\int_0^1|S_1(t)|^{2l}{\rm d}t=\mcal{O}_l(1),\ T\to\infty, 
\ee  
then it is true that\footnote{See (\ref{3.1}) and (\ref{3.11}).} 
\be \label{3.18} 
\int_1^T\left\{\frac{1}{\zeta(2\sigma)}|\zeta(\sigma+it)|^2+\frac{1}{\bar{c}(l)}|S_1(t)|^{2l}\right\}{\rm d}\sim 2T, 
\ee  
i. e. 
\be \label{3.19} 
\int_1^T\left\{ \bar{c}(l)|\zeta(\sigma+it)|^2+\zeta(2\sigma)|S_1(t)|^{2l}\right\}{\rm d}t\sim 2\bar{c}(l)\zeta(2\sigma)T. 
\ee  
Now, if we use the substitution 
\be \label{3.20} 
T=\frac{x}{2\bar{c}(l)\zeta(2\sigma)}\tau,\ x>0, 
\ee  
in (\ref{3.19}), then we obtain the following functional. 

\begin{mydef55}
\be \label{3.21} 
\begin{split}
& \lim_{\tau\to\infty}\frac{1}{\tau}\int_1^{\frac{x}{2\bar{c}(l)\zeta(2\sigma)}\tau}\{\bar{c}(l)|\zeta(\sigma+it)|^2+\zeta(2\sigma)|S_1(t)|^2\}{\rm d}t=x 
\end{split}
\ee 
for every fixed $x>0$, $\sigma>1$ and $l\in\mbb{N}$. 
\end{mydef55} 

Next, if we use the substitution (\ref{3.7}) in (\ref{3.21}), then we obtain the following. 

\begin{mydef56}
\be \label{3.22} 
\begin{split}
& \lim_{\tau\to\infty}\frac{1}{\tau}\int_1^{\frac{x^n+y^n}{z^n}\frac{\tau}{2\bar{c}(l)\zeta(2\sigma)}}
\{\bar{c}(l)|\zeta(\sigma+it)|^2+\zeta(2\sigma)|S_1(t)|^2\}{\rm d}t= \\ 
& \frac{x^n+y^n}{z^n} 
\end{split}
\ee 
for every fixed Fermat's rational.  
\end{mydef56} 

Consequently, we get the result. 

\begin{mydef13}
The $\zeta$-condition 
\be \label{3.23} 
\begin{split}
	& \lim_{\tau\to\infty}\frac{1}{\tau}\int_1^{\frac{x^n+y^n}{z^n}\frac{\tau}{2\bar{c}(l)\zeta(2\sigma)}}
	\{\bar{c}(l)|\zeta(\sigma+it)|^2+\zeta(2\sigma)|S_1(t)|^2\}{\rm d}t\not= 1 
\end{split}
\ee 
on the class of all Fermat's rationals represents the next $\zeta$-equivalent of the Fermat-Wiles theorem for every fixed $\sigma>1$ and $l\in\mbb{N}$. 
\end{mydef13} 

\begin{remark}
The $\zeta$-conditions (\ref{3.9}) and (\ref{3.23}) can also be generalized as it follows from the Remark 5. 
\end{remark} 

\section{Jacob's ladders and the second set of equivalents of the Fermat-Wiles theorem} 

\subsection{} 

In this section we use, firstly, the formula\footnote{See (\ref{1.1}) and (\ref{1.3}), comp. \cite{16}, pp. 30, 31.}
\bdis  
\int_1^T|\zeta(\sigma+it)|^2{\rm d}t=\zeta(2\sigma)T+R_a(T,\sigma),\ \sigma\geq \frac 12+\epsilon, 
\edis  
where 
\bdis 
R_1(T,\sigma)=\mcal{O}(T^{2-2\sigma}\ln T)|_{\sigma\geq \frac{1}{2}+\epsilon}=\mcal{O}(T^{1-2\epsilon}\ln T) 
\edis  
i. e. 
\be \label{4.1} 
\int_1^T|\zeta(\sigma+it)|^2{\rm d}t=\zeta(2\sigma)T+\mcal{O}(T^{1-2\epsilon}\ln T),\ \sigma\geq\frac 12+\epsilon, 
\ee  
for sufficiently small positive $\epsilon$.  

\begin{remark}
The essential part of the formula (\ref{4.1}), that is the part corresponding to the interval 
\bdis 
\frac 12<\sigma\leq 1, 
\edis  
has been obtained by Hardy and Littlewood. 
\end{remark} 

Next, we use the first Selberg's formula, see (\ref{3.10}) 
\be \label{4.2} 
\begin{split}
& \int_T^{T+H}|S_1(t)|^{2l}{\rm d}t=\bar{c}(l)H+\mcal{O}\left(\frac{H}{\ln T}\right),\ T^a\leq H\leq T,\ \frac 12<a\leq 1. 
\end{split} 
\ee 
And finally, we use our almost linear formula 
\be \label{4.3} 
\int_{\overset{r-1}{T}}^{\overset{r}{T}}\left|\zf\right|^2{\rm d}t=(1-c)\overset{r-1}{T}+\mcal{O}(T^{\frac 13+\delta}),\ r=1,\dots,k, 
\ee  
for every fixed\footnote{See \cite{8}, (3.4) and (3.6).} $k\in\mbb{N}$ and $0<\delta$ being sufficiently small.  

\subsection{} 

Next, the reversal Jacob's ladders generate the following increments: 
\begin{itemize}
\item[(A)] for the formula (\ref{4.1}) we have\footnote{See (\ref{2.5}), (\ref{2.14}).} 
\be \label{4.4} 
\int_{\overset{r-1}{T}}^{\overset{r}{T}}|\zeta(\sigma+it)|^2{\rm d}t=\zeta(2\sigma)(\overset{r}{T}-\overset{r-1}{T})+\mcal{O}(T^{1-2\epsilon}\ln T),\ r=1,\dots k; 
\ee 	
\item[(B)] for the formula (\ref{4.2}) we have\footnote{See again (\ref{2.5}), (\ref{2.14}).} 
\be \label{4.5} 
\int_{\overset{r-1}{T}}^{\overset{r}{T}}|S_1(t)|^{2l}{\rm d}t=\bar{c}(l)(\overset{r}{T}-\overset{r-1}{T})+\mcal{O}\left(\frac{T}{\ln^2T}\right), 
\ee  
here, of course, 
\bdis 
\overset{r}{T}=\overset{r-1}{T}+(\overset{r}{T}-\overset{r-1}{T})=\overset{r-1}{T}+H. 
\edis 
\end{itemize} 

Now, we have the following result\footnote{See (\ref{4.3}) -- (\ref{4.5}).}. 

\begin{mydef57}
\be \label{4.6} 
\int_{\overset{r-1}{T}}^{\overset{r}{T}}\left|\zf\right|^2{\rm d}t\sim (1-c)\overset{r-1}{T}, 
\ee  
\be \label{4.7} 
\int_{\overset{r-1}{T}}^{\overset{r}{T}}|\zeta(\sigma+it)|^2{\rm d}t\sim\zeta(2\sigma)(\overset{r}{T}-\overset{r-1}{T}), 
\ee 
\be \label{4.8} 
\int_{\overset{r-1}{T}}^{\overset{r}{T}}|S_1(t)|^{2l}{\rm d}t\sim \bar{c}(l)(\overset{r}{T}-\overset{r-1}{T}), 
\ee 
\bdis 
r=1,\dots,k,\ \overset{r}{T}=\overset{r}{T}(T)=\vp_1^{-r}(T),\ T\to\infty 
\edis 
for every fixed $k,l\in\mbb{N}$, $\sigma\geq \frac 12+\epsilon$. 
\end{mydef57}  

\subsection{} 

It follows from (\ref{4.6}) and (\ref{4.7}) that 
\bdis 
\int_{\overset{r-1}{T}}^{\overset{r}{T}}\left\{ \frac{1}{1-c}\left|\zf\right|^2+\frac{1}{\zeta(2\sigma)}|\zeta(\sigma+it)|^2\right\}{\rm d}t\sim \overset{r}{T}, 
\edis 
i. e. 
\be \label{4.9} 
\begin{split}
& \int_{\overset{r-1}{T}}^{\overset{r}{T}}\left\{ \zeta(2\sigma)\left|\zf\right|^2+(1-c)|\zeta(\sigma+it)|^2\right\}{\rm d}t\sim \\ 
& (1-c)\zeta(2\sigma)\overset{r}{T},\ r=1,\dots,k,\ T\to\infty. 
\end{split}
\ee  
By summation in (\ref{4.9}) we obtain 
\be \label{4.10} 
\begin{split}
	& \int_{T}^{\overset{k}{T}}\left\{ \zeta(2\sigma)\left|\zf\right|^2+(1-c)|\zeta(\sigma+it)|^2\right\}{\rm d}t\sim \\ 
	& (1-c)\zeta(2\sigma)\left(\overset{1}{T}+\overset{2}{T}+\dots+\overset{k}{T}\right). 
\end{split}
\ee 

Now, it follows\footnote{See (\ref{2.5}) and (\ref{2.14}).} that
\be \label{4.11} 
\overset{1}{T}+\overset{2}{T}+\dots+\overset{k}{T}=kT+\mcal{O}\left( k^2\frac{T}{\ln T}\right),\ T\to\infty
\ee 
for every fixed $k\in\mbb{N}$. Consequently, we have the following result. 

\begin{mydef58}
\be \label{4.12} 
\begin{split}
& \int_{T}^{\overset{k}{T}}\left\{ \zeta(2\sigma)\left|\zf\right|^2+(1-c)|\zeta(\sigma+it)|^2\right\}{\rm d}t\sim \\ 
& k(1-c)\zeta(2\sigma)T,\ T\to\infty 
\end{split}
\ee 
for every fixed $k\in\mbb{N}$ and $\sigma\geq \frac 12+\epsilon$. 
\end{mydef58} 

Next, the substitution 
\be \label{4.13} 
T=\frac{x}{k(1-c)\zeta(2\sigma)}\tau,\ x>0 
\ee  
in (\ref{4.12}) gives the following functional. 

\begin{mydef59}
\be \label{4.14}
\begin{split}
& \lim_{\tau\to\infty}\frac{1}{\tau}\int_{\frac{x}{k(1-c)\zeta(2\sigma)}\tau}^{[\frac{x}{k(1-c)\zeta(2\sigma)}\tau]^k}\left\{ \zeta(2\sigma)\left|\zf\right|^2+(1-c)|\zeta(\sigma+it)|^2\right\}{\rm d}t=x
\end{split}
\ee 
for every fixed $x>0$, $k\in\mbb{N}$ and $\sigma\geq \frac 12+\epsilon$, where 
\bdis 
[Y]^k=\vp_1^{-k}(Y). 
\edis 
\end{mydef59} 

Now, if we use the substitution (\ref{3.7}) in (\ref{4.14}), then we obtain the following 
\begin{mydef510}
\be \label{4.15} 
\begin{split}
& \lim_{\tau\to\infty}\frac{1}{\tau}\int_{\frac{x^n+y^n}{z^n}\frac{\tau}{k(1-c)\zeta(2\sigma)}}^{[\frac{x^n+y^n}{z^n}\frac{\tau}{k(1-c)\zeta(2\sigma)}]^k}\left\{ \zeta(2\sigma)\left|\zf\right|^2+(1-c)|\zeta(\sigma+it)|^2\right\}{\rm d}t \\ 
& = \frac{x^n+y^n}{z^n}
\end{split}
\ee 
for every fixed Fermat's rational. 
\end{mydef510} 

Consequently, we have the next result. 

\begin{mydef14}
The $\zeta$-condition 
\be \label{4.16} 
\begin{split}
& \lim_{\tau\to\infty}\frac{1}{\tau}\int_{\frac{x^n+y^n}{z^n}\frac{\tau}{k(1-c)\zeta(2\sigma)}}^{[\frac{x^n+y^n}{z^n}\frac{\tau}{k(1-c)\zeta(2\sigma)}]^k}\left\{ \zeta(2\sigma)\left|\zf\right|^2+(1-c)|\zeta(\sigma+it)|^2\right\}{\rm d}t \\ 
& \not=1
\end{split}
\ee 
on the class of all Fermat's rationals represents the next $\zeta$-equivalent of the Fermat-Wiles theorem for every fixed $k\in\mbb{N}$ and $\sigma\geq \frac 12+\epsilon$. 
\end{mydef14} 

\subsection{} 

In the case of formulae (\ref{4.6}) and (\ref{4.8}) we use the method presented in the subsection 4.3 and we obtain the following. 

\begin{mydef15}
The $\zeta$-condition 
\be \label{4.17} 
\begin{split}
	& \lim_{\tau\to\infty}\frac{1}{\tau}\int_{\frac{x^n+y^n}{z^n}\frac{\tau}{k(1-c)\bar{c}(l)}}^{[\frac{x^n+y^n}{z^n}\frac{\tau}{k(1-c)\bar{c}(l)}]^k}\left\{ \bar{c}(l)\left|\zf\right|^2+(1-c)|S_1(t)|^{2l}\right\}{\rm d}t \\ 
	& \not=1
\end{split}
\ee 
on the class of all Fermat's rationals represents the next $\zeta$-equivalent of the Fermat-Wiles theorem for every fixed $k,l\in\mbb{N}$. 
\end{mydef15} 

\subsection{} 

It follows from (\ref{4.6}) -- (\ref{4.8}) that 
\bdis 
\begin{split}
& \int_{\overset{r-1}{T}}^{\overset{r}{T}}
\left\{ 
\frac{2}{1-c}\left|\zf\right|^2+\frac{1}{\zeta(2\sigma)}|\zeta(\sigma+it)|^2+\frac{1}{\bar{c}(l)}|S_1(t)|^{2l}
\right\}{\rm d}t\sim 2\overset{r}{T}, 
\end{split}
\edis  
i. e. 
\be \label{4.18} 
\begin{split}
& \int_{\overset{r-1}{T}}^{\overset{r}{T}}{\rm d}t\times \\ 
& \left\{
2\bar{c}(l)\zeta(2\sigma)\left|\zf\right|^2+(1-c)\bar{c}(l)|\zeta(\sigma+it)|^2+(1-c)\zeta(2\sigma)|S_1(t)|^{2l}
\right\}\sim \\ 
& 2(1-c)\bar{c}(l)\zeta(2\sigma)\overset{r}{T}, \\ 
& r=1,\dots,k,\ T\to\infty. 
\end{split}
\ee 
Next, it follows from (\ref{4.18}), comp. (\ref{4.10}) and (\ref{4.11}). 

\begin{mydef511}
\be \label{4.19} 
\begin{split}
& \int_T^{\overset{k}{T}}
\left\{ 
2\bar{c}(l)\zeta(2\sigma)\left|\zf\right|^2+(1-c)\bar{c}(l)|\zeta(\sigma+it)|^2+ \right. \\ 
& \left. (1-c)\zeta(2\sigma)|S_1(t)|^2
\right\}{\rm d}t\sim 2k(1-c)\bar{c}(l)\zeta(2\sigma)T,\ T\to\infty. 
\end{split}
\ee 
\end{mydef511} 

Consequently, the method (\ref{4.13}) -- (\ref{4.16}) applied on the formula (\ref{4.19}) gives the following result. 

\begin{mydef16}
The $\zeta$-condition 
\be \label{4.20} 
\begin{split}
& \lim_{\tau\to\infty}\frac{1}{\tau}\int_{\frac{x^n+y^n}{z^n}\frac{\tau}{2k(1-c)\bar{c}(l)\zeta(2\sigma)}}^{[\frac{x^n+y^n}{z^n}\frac{\tau}{2k(1-c)\bar{c}(l)\zeta(2\sigma)}]^k}\left\{ 
2\bar{c}(l)\zeta(2\sigma)\left|\zf\right|^2+ \right. \\ 
& \left. (1-c)\bar{c}(l)|\zeta(\sigma+it)|^2+  (1-c)\zeta(2\sigma)|S_1(t)|^2
\right\}{\rm d}t\not= 1
\end{split}
\ee 
on the class of all Fermat's rationals represents the next $\zeta$-equivalent of the Fermat-Wiles theorem for every fixed 
\bdis 
k,l\in\mbb{N},\ \sigma\geq \frac{1}{2}+\epsilon. 
\edis 
\end{mydef16}

\begin{remark}
Of course, the $zeta$-conditions (\ref{4.16}) and (\ref{4.20}) can be generalized in the direction of Remark 5. 
\end{remark}

\section{On the existence of infinite set of point of contact between the set of all Dirichlet's series and the Fermat-Wiles theorem} 

\subsection{} 

Following our subsection 1.5 we have by eqs. (\ref{1.24}) and (\ref{1.25}) the following 
\be \label{5.1} 
\int_0^T|f(\sigma_0+it)|^2{\rm d}t=F(\sigma_0;f)T+o(T),\ T\to\infty. 
\ee  
Now, if we put 
\be \label{5.2} 
T=\frac{x}{F(\sigma_0;f)}\tau,\ x>0
\ee 
into (\ref{5.1}), then we obtain this result. 

\begin{mydef512}
\be \label{5.3} 
\lim_{\tau\to\infty}\frac{1}{\tau}\int_0^{\frac{x}{F(\sigma_0;f)}\tau}|f(\sigma_0+it)|^2=x
\ee 
for every fixed $x>0$, $f(s)\in\mfrak{D}$ and $\sigma_0(f)$. 
\end{mydef512} 

Next, by the substitution (\ref{3.7}) in (\ref{5.3}), we instantly obtain the following statement. 

\begin{mydef513}
\be \label{5.4} 
\lim_{\tau\to\infty}\frac{1}{\tau}\int_0^{\frac{x^n+y^n}{z^n}\frac{\tau}{F(\sigma_0;f)}}|f(\sigma_0+it)|^2{\rm d}t=\frac{x^n+y^n}{z^n}
\ee  
for every fixed Fermat's rational and every fixed 
\bdis 
f(s)\in\mfrak{D},\ \sigma_0(f). 
\edis 
\end{mydef513} 

Consequently, we obtain the following theorem. 

\begin{mydef17} 
The $\mfrak{D}$-condition 
\be \label{5.5} 
\lim_{\tau\to\infty}\frac{1}{\tau}\int_0^{\frac{x^n+y^n}{z^n}\frac{\tau}{F(\sigma_0;f)}}|f(\sigma_0+it)|^2{\rm d}t\not=1 
\ee  
on the class of all Fermat's rationals represents new $\mfrak{D}$-equivalent of the Fermat-Wiles theorem for every fixed $f(s)\in\mfrak{D}$ and $\sigma_0(f)$. 
\end{mydef17} 

I would like to thank Michal Demetrian for his moral support of my study of Jacob's ladders.

\end{document}